\documentclass[utf8]{FrontiersinHarvard}
\usepackage[utf8]{inputenc}
\usepackage[T1]{fontenc}
\usepackage{graphicx}
\usepackage{longtable}
\usepackage{wrapfig}
\usepackage{rotating}
\usepackage[normalem]{ulem}
\usepackage{amsmath}
\usepackage{amssymb}
\usepackage{capt-of}
\usepackage{hyperref}
\usepackage{amsmath}                                  
\usepackage{amssymb}                                  
\usepackage{booktabs}                                 
\usepackage{cite}                                     
\usepackage{graphicx}                                 
\usepackage{lineno}                                   
\usepackage{microtype}                                
\usepackage{multicol}                                 
\usepackage{multirow}                                 
\usepackage{pgfplots}                                 
\usepackage{subcaption}                               
\usepackage{subfloat}                                 
\usepackage{soul}                                     
\usepackage{tabularx}                                 
\usepackage{tikz}                                     
\usepackage{xcolor}                                   
\usepackage{xfp}                                      
\usepackage{standalone}                               
\usepackage{hyperref}                                 
\usetikzlibrary{automata, positioning, arrows.meta}   
\pgfplotsset{compat=1.3, width=\textwidth}
\graphicspath{ {img} }                                
\let\ref\autoref                                      
\let\cite\citep                                       

\newcommand{\A}{35 }                                                            
\newcommand{\N}{338 }                                                           
\newcommand{\acharge}{0.9}                                                      
\newcommand{\bcharge}{0.7 }                                                     
\newcommand{\mincharge}{20\% }                                                  
\newcommand{\minchargeD}{0.20 }                                                 
\newcommand{\batsize}{388 }                                                     
\newcommand{\fast}{15 }                                                         
\newcommand{\slow}{15 }                                                         
\newcommand{\fasts}{911 }                                                       
\newcommand{\slows}{30 }                                                        
\newcommand{\contvars}{7,511 }
\newcommand{\intvars}{328,282 }
\newcommand{\timeran}{7200 }                                                    
\def\keyFont{\fontsize{8}{11}\helveticabold }
\def\firstAuthorLast{Brown {et~al.}} 
\def\Authors{Alexander Brown\,$^{1,*}$, Greg Droge\,$^{1}$, Jacob Gunther\,$^{1}$}
\def\Address{$^{1}$Department of Electrical and Computer Engineering, Logan, UT, USA}
\def\corrAuthor{Alexander Brown}
\def\corrEmail{A01704744@usu.edu}
\author[\firstAuthorLast ]{\Authors} 
\address{} 
\correspondance{} 
\extraAuth{}
\date{}
\title{}
\hypersetup{
 pdfauthor={},
 pdftitle={},
 pdfkeywords={},
 pdfsubject={},
 pdfcreator={Emacs 29.3 (Org mode 9.6.15)},
 pdflang={English}}
\begin{document}

\title{A Position Allocation Approach to the Scheduling of Battery-Electric Bus Charging}
\onecolumn
\firstpage{1}

\maketitle

\begin{abstract}
  Robust charging schedules in a growing market of battery electric bus (BEB) fleets are a critical component to
  successful adoption. In this paper, a BEB charging scheduling framework that considers spatiotemporal schedule
  constraints, route schedules, fast and slow charging, and battery dynamics is modeled as a mixed integer linear
  program (MILP). The MILP is modeled after the Berth Allocation Problem (BAP) in a modified form known as the Position
  Allocation Problem (PAP). Linear battery dynamics are included to model the charging of buses while at the station. To
  model the BEB discharges over their respective routes, it is assumed each BEB has an average kWh charge loss while on
  route. The optimization coordinates BEB charging to ensure that each vehicle remains above a specified state-of-charge
  (SOC). The model also minimizes the total number of chargers utilized and prioritizes slow charging for battery
  health. The model validity is demonstrated with a set of routes sampled from the Utah Transit Authority (UTA) for \A
  buses and \N visits to the charging station. The model is also compared to a heuristic algorithm based on charge
  thresholds referred to as the Qin-Modified method. The results presented show that the slow chargers are more readily
  selected and the charging and spatiotemporal constraints are met while considering the battery dynamics and minimizing
  both the charger count and consumption cost.

  \tiny \keyFont{ \section{Keywords:} Berth Allocation Problem (BAP), Position Allocation Problem (PAP), Mixed Integer
    Linear Program (MILP), Battery Electric Bus (BEB), Scheduling}
\end{abstract}
\section{Introduction}
\label{sec:introduction}
The public transportation system is crucial in any urban area; however, the increased awareness and concern of the
environmental impacts of petroleum-based public transportation has driven an effort to reduce the pollutant footprint
\cite{de-2014-simul-elect,xylia-2018-role-charg,guida-2017-zeeus-repor-europ,li-2016-batter-elect}. Particularly,
the electrification of public bus transportation via battery power, i.e., battery-electric buses (BEBs), have received
significant attention \cite{li-2016-batter-elect}. Although the technology provides benefits beyond a reduction in
emissions, such as lower driving costs, lower maintenance costs, and reduced vehicle noise, battery-powered systems
introduce new challenges such as larger upfront costs, and potentially several hours long ``refueling'' periods
\cite{xylia-2018-role-charg,li-2016-batter-elect}. Furthermore, the problem is exacerbated by the constraints of the
transit schedule to which the fleet must adhere, the limited amount of chargers available, and the adverse effects on
the health of the battery due to fast charging \cite{lutsey-2019-updat-elect}. This paper presents a
framework for optimally assiging BEBs to charging queues assuming fixed routes while taking into consideration multiple
charger types and utilizing linear charging dynamics. This method also enforces the SOC to stay above a specified
percentage throughout the day, and ensures a minumum SOC at the end of the working day.

Many recent efforts have been made to simultaneously solve the problems of route scheduling, charging fleets, and
determining the infrastructure upon which they rely, e.g., \cite{wei-2018-optim-spatio,sebastiani-2016-evaluat-elect,hoke-2014-accoun-lithium,wang-2017-elect-vehic}. Several simplifications are made to make these problems
computationally feasible. Simplifications to the charge scheduling model include utilizing only fast chargers while
planning \cite{wei-2018-optim-spatio,sebastiani-2016-evaluat-elect,wang-2017-optim-rechar,zhou-2020-bi-objec,yang-2018-charg-sched,wang-2017-elect-vehic,qin-2016-numer-analy,liu-2020-batter-elect}. If slow chargers are used,
they are only employed at the depot and not the station \cite{he-2020-optim-charg,tang-2019-robus-sched}. Some
approaches also simplify by assuming a full charge is always achieved
\cite{wei-2018-optim-spatio,wang-2017-elect-vehic,zhou-2020-bi-objec,wang-2017-optim-rechar}. Others have assumed
that the charge received is proportional to the time spent on the charger
\cite{liu-2020-batter-elect,yang-2018-charg-sched}, which can be a valid assumption when the battery state-of-charge
(SOC) is below 80\% \cite{liu-2020-batter-elect}.

This work builds upon the Position Allocation Problem (PAP) \cite{qarebagh-2019-optim-sched}, a modification of the
well-studied Berth Allocation Problem (BAP), as a means to schedule the charging of electric vehicles
\cite{buhrkal-2011-model-discr,frojan-2015-contin-berth,imai-2001-dynam-berth}. The BAP is a continuous time model
that solves the problem of allocating space for incoming vessels to be berthed and serviced. Each arriving vessel
requires both time and space to be serviced and thus must be carefully assigned to avoid delay
\cite{imai-2001-dynam-berth}. Vessels are lined up parallel to the berth to be serviced and are horizontally queued as
shown in \autoref{subfig:bapexample}. As the vessels are serviced, they move from left to right to make space for the
queued vessels moving vertically downward into their respective berthing locations. The PAP utilizes this notion of
queuing for scheduling vehicles to be charged, as shown in \autoref{subfig:papexample}. The vehicles are queued in
several lines and move from left to right to recieve their charge and exit the system. The PAP is formulated as a
rectangle packing problem and assumes that each vehicle has a predefined charge time, the amount of vehicles that can
charge at any given moment is limited by the physical width of each vehicle and the length of the charging block. The
PAP also makes the assumption that each vehicle that is placed in the system is unique
\cite{qarebagh-2019-optim-sched}.

The main contribution of this work is the extension of the PAP's novel approach to BEB charger scheduling. This
incorporates a proportional charging model into the MILP framework, includes consideration for multiple charger types,
and consideration of each route in the schedule. The last contribution is of importance because both the BAP and PAP
consider each arrival to be unique; thus, the tracking of battery charge from one visit to the next must be considered.
Furthermore, the input parameters for the model can be predefined in such a manner as to minimize the number of fast and
slow chargers utilized as well as minimize the energy consumption. That is, the model will simultaneously minimize the
number of chargers as well as the total consumed energy. The result is a MILP formulation that coordinates charging
times and charger type for every visit while considering a dynamic charge model with scheduling constraints.

The remainder of the paper proceeds as follows: In \autoref{sec:the-position-allocation-problem}, the PAP is introduced
with a formulation of the resulting MILP. \autoref{sec:problemformulation} constructs the MILP for BEB scheduling,
including modifications to the PAP queuing constraints and the development of a dynamic charging model.
\autoref{sec:example} demonstrates an example of using the formulation to coordinate \A buses over \N total visits to
the station. The paper ends in \autoref{sec:conclusion} with concluding remarks.
\section{The Position Allocation Problem}
\label{sec:the-position-allocation-problem}
This section provides a brief overview of the BAP and a detailed formulation of PAP as presented in
\cite{qarebagh-2019-optim-sched}.

\subsection{Overview of BAP}
\label{sec:overview-of-bap}
\begin{figure}[htpb]
\centering
    \includegraphics{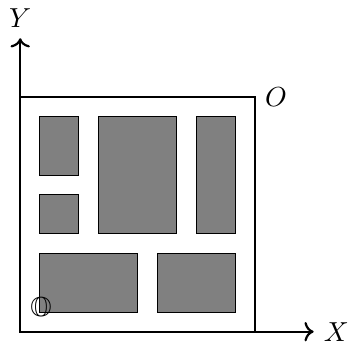}
    \caption{Example of the rectangle packing problem. The large square represented by $O$ indicates the constrained
      area that the set of shaded rectangles $\mathbb{O}$ must be placed within.}
    \label{fig:packexample}
\end{figure}

The BAP is a rectangle packing problem where a set of rectangles, \(\mathbb{O}\), are attempted to be optimally placed in
a larger rectangle, \(O\), as shown in \autoref{fig:packexample}. The rectangle packing problem is an NP-hard problem that
can be used to describe many real-life problems \cite{bruin-2013-rectan-packin,murata-1995-rectan}. In some of these
problems, the dimensions of \(\mathbb{O}\) are held constant such as in the problem of packing modules on a chip, where
the widths and height of the rectangles represent the physical width and heights of the modules
\cite{murata-1995-rectan}. Other problems, such as the one presented in this work, allow either the horizontal or
vertical edge of each rectangle in \(\mathbb{O}\) to vary. As an example, suppose the vessel lengths are predefined
(vertical edges are static), but the service time is allowed to vary (horizontal edges are dynamic).
\cite{buhrkal-2011-model-discr}.

\begin{figure}[ht]
\centering
    \includegraphics{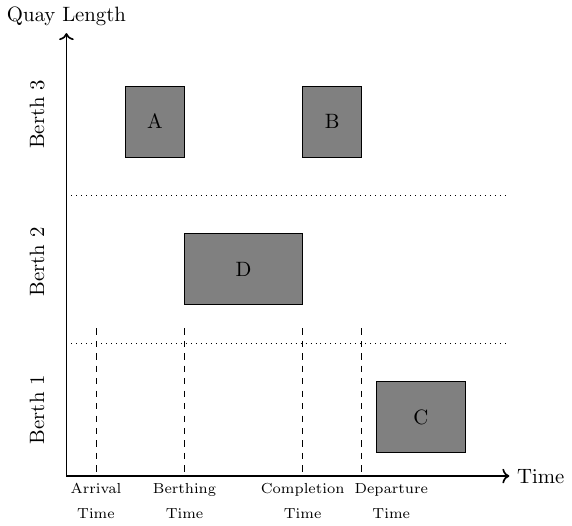}
    \caption{The representation of the berth-time space. The x and y-axis represent time and space, respectively. Along
      the y-axis, the dashed lines represent discrete berthing locations. These locations may be chosen to be
      continuous. The shaded rectangles represent scheduled vessels to be serviced. The height of each shaded rectangle
      represents the space taken on the berth and the width being the time to service said vessel. The vertical dashed
      lines are associated with vessel D and represent the arrival time, berthing time, service completion time, and
      departure time. Note that the arrival time may be before the berthing time and the completion time may be before
      the departure time.}
    \label{fig:bap}
\end{figure}

The BAP solves the problem of optimally assigning incoming vessels to berth positions to be serviced as shown in
\autoref{subfig:bapexample}. To relate to the rectangle packing problem, the width and height of \(O\) represent the time
horizon \(T\) seconds and the berth length \(L\) meters, respectively. Similarly, the widths and heights of each element in
\(\mathbb{O}\) represent the time spent to service vessel \(i\) and the space taken by docking vessel \(i\), respectively. In
the BAP, the vessel characteristics (length of the vessel, arrival time, handling time, desired departure time) are
assumed to be known for all vessels. A representation of a BAP solution is shown in \autoref{fig:bap}. The x and y-axis
represent time horizon and berthing space, respectively. The gray squares, labeled A, B, C, and D, represent berthed
vessels. The width of the boxes represents the time spent being serviced, and the height represents the amount of space
the vessel requires on the berth. The vertical line adjacent to ``Arrival Time'' represents the actual time that the
vessel arrives and is available to be berthed. ``Berthing Time'' is the time the vessel is berthed and begins being
serviced. ``Completion time'' represents the time at which the berthing space becomes available again.

\subsection{The PAP Formulation}
\label{sec:the-pap-formulation}
The BAP forms the basis of the PAP; however, there are some differences in the way the variables are interpreted. The
starting service time, \(u_i\) seconds, is viewed as the initial charge time, and the service time, total elapsed time
spent on the charger. Similarly, for the spatial term, \(v_i \in [0,L]\), the berth location is instead interpreted as the
initial position on the charger. There are also a few clarifying concepts about how the system is modeled. The PAP
models the set of chargers as one continuous line; that is, the natural behavior of the PAP model is to allow vehicles
to be queued anywhere along \([0,L]\). Similarly, the charge times are continuous and can be placed anywhere on the time
horizon, \([0,T]\), as long as the allocated times do not interfere with other scheduled charge times.

The PAP formulation's parameters can be divided into two categories: input parameters and decision variables. Each type
will now be introduced in turn. The following parameters are assumed to be known inputs for the MILP. \(L\) defines the
length of the charger in meters. As stated previously, it is modeled as a continuous bar meaninng that a vehicle can be
placed anywhere in the range \([0,L]\). It is assumed that the time horizon, \(T\) seconds, is known so that vehicles may be
placed temporarily in the range \([0,T]\). The total number of visits to the station over the time horizon is represented
by \(n_V\). The arrival time for each visit is represented by \(a_i\) seconds, and the required charge time is represented
by \(s_i\) seconds. The width of vehicle \(i\) is represented by \(l_i\) meters.

The decision variables provide the means by which the solver may optimize the problem. The initial and final charge
times for vehicle \(i\) are \(u_i\) and \(d_i\) seconds, respectively. The starting position on the charger is denoted as \(v_i
\in [0,L]\) meters. The temporal ordering of vehicles \(i\) and \(j\) is determined by \(\sigma_{ij} \in \{0, 1\}\), where \(\sigma_{ij} = 1
\implies\) \(i\) arrives before \(j\) for all \(1 \le i,j \le n_V\). Similarly, \(\psi_{ij} \in \{0, 1\}\) determines the relative
position of vehicles \(i\) and \(j\) on the charger: \(\psi_{ij} = 1 \implies v_i < v_j\) for all \(1 \le i,j \le n_V\).

To determine the values for each of these decision variables, a MILP was formulated in
\cite{qarebagh-2019-optim-sched}. The formulation is shown in its entirety for completeness.
The problem to be solved is

\begin{equation}
	\label{eq:bapobjective}
	\min\; \sum_{i=1}^N (d_i - a_i)
\end{equation}

Subject to:
\begin{subequations}
\label{eq:bapconstrs}
\begin{align}
    u_j - u_i - s_i - (\sigma_{ij} - 1)T \geq 0   \label{subeq:baptime}          \\
    v_j - v_i - l_i - (\psi_{ij} - 1)L \geq 0   \label{subeq:bapspace}           \\
    \sigma_{ij} + \sigma_{ji} + \psi_{ij} + \psi_{ji} \geq 1 \label{subeq:bapvalid_pos}     \\
    \sigma_{ij} + \sigma_{ji} \leq 1                   \label{subeq:bapsigma}        \\
    \psi_{ij} + \psi_{ji} \leq 1                   \label{subeq:bapdelta}        \\
    s_i + u_i = d_i                       \label{subeq:bapdetach}       \\
    a_i \leq u_i \leq (T - s_i)                 \label{subeq:bapvalid_starts} \\
    \sigma_{ij} \in \{0,1\},\;\psi_{ij} \in \{0,1\}\; \label{subeq:bapsdspace}      \\
    v_i \in [0, L ]                         \label{subeq:bapvspace}
\end{align}
\end{subequations}

\noindent The objective function, \autoref{eq:bapobjective}, minimizes the idle and service time by summing over the
differences between the departure time, \(d_i\), and arrival time, \(a_i\) for all visits. In other words, the objective
function is searching for the schedule that removes each vehicle from the service queue as quickly as possible.

\autoref{subeq:baptime}-\autoref{subeq:bapdelta} are used to ensure that individual rectangles do not overlap. In terms
of the PAP, this implies that there are no conflicts in the schedule spatially or temporally. \autoref{subeq:baptime}
establishes temporal ordering when active (\(\sigma_{ij}=1\)) in the manner described previously by utilizing big-M notation.
Similarly, \autoref{subeq:bapspace} establishes spatial ordering when active (\(\psi_{ij} =1\)). Constraints
\autoref{subeq:bapvalid_pos}-\autoref{subeq:bapdelta} enforce spatial and temporal ordering between each queue/vehicle
pair. Constraints \autoref{subeq:bapsigma} and \autoref{subeq:bapdelta} enforce validity of the assignments. For
example, if \autoref{subeq:bapsigma} resulted in a value of two, that would imply both vehicle \(i\) and \(j\) are scheduled
before and after each other temporally, which is impossible. In the case of \autoref{subeq:bapdelta} being equal to
two, that would mean that vehicles \(i\) and \(j\) are scheduled both before and after one another on the charging strip,
which is again impossible.

The last constraints force relationships between arrival time, initial charge time, and departure time.
\autoref{subeq:bapdetach} states that the initial charge time, \(u_i\), plus the total charge time for, \(s_i\), must equal
the departure time, \(d_i\). \autoref{subeq:bapvalid_starts} enforces the arrival time, \(a_i\), to be less than or equal to
the service start time, \(u_i\), which in turn must be less than or equal to the latest time the vehicle may begin
charging and stay within the time horizon. \autoref{subeq:bapsdspace} simply states that \(\sigma_{ij}\) and \(\psi_{ij}\) are
binary terms. \autoref{subeq:bapvspace} ensures that the assigned value of \(v_i\) is within the range, \([0,L]\).
\section{A Rectangle Packing Formulation for BEB Charging}
\label{sec:problemformulation}
Applying the PAP to BEB charging requires four fundamental changes. The first is that the time that a BEB spends
charging must be allowed to vary. That is, \(u_i\), \(d_i\), and \(s_i\) become variables of optimization. This is done
primarily because chargers of various speeds are to be introduced. Allowing BEBs to have multiple visits that a charger
decided upon during the optimization requires that the start and stop times must be changeable to respect the SOC
constraints. Second, in the PAP each visit is assumed to be a different vehicle. For the BEB charging problem, each bus
may make multiple visits to the station throughout the day. Thus, the resulting SOC for a bus at a given visit is
dependent upon each of the prior visits. The third fundamental change is related to the first two. The SOC of each bus
must be tracked to ensure that charging across multiple visits is sufficient to allow each bus to execute its route
throughout the day. Finally, as previously stated, the PAP models the charger as one continuous bar. For the BEB, it
will be assumed that a discrete number of chargers exist. Moreover, it is assumed that these chargers may have different
charge rates.

A few assumptions are made in the derivation of the algorithm. As this work is not focused on estimating the discharge
of a BEB during its route, the discharge for each route will be pre-calculated by assuming a fixed discharge rate kW
multiplied by the route duration in hours. Secondly, it is assumed that the initial SOC of each BEB at the beginning of
the day, $\alpha_b\kappa_b$, is larger than the minimum required SOC at the end of the day, $\beta_b\kappa_b$.
Therefore, it must be assumed that the difference in the SOC can reach $\alpha_b\kappa_b$ by the beginning of the next
working day.

The discussion of the four changes is separated into two sections. \autoref{sec:queuing} discusses the changes in the
spatial-temporal constraint formulation to form a queuing constraint. \autoref{sec:batt_dynamics} then discusses the
addition of bus charge management. This section ends with a brief discussion of a modified objective function and the
statement of the full problem in \autoref{sec:BEB_MILP}. The notation is explained throughout and summarized in
\autoref{tab:variables}.

\subsection{Queuing Constraints}
\label{sec:queuing}
\noindent The queuing constraints ensure that the buses entering the charging queues are assigned
feasibly. There are three sets to differentiate between different entities. \(\mathbb{B} = \{1, ..., n_B\}\) is the set of
bus indices with index \(b\) used to denote an individual bus, \(\mathbb{Q} = \{1, ..., n_Q\}\) is the set of queues with index \(q\)
used to denote an individual queue, and \(\mathbb{V} = \{1, ..., n_V\}\) is a set of visits to the station with \(i\) and
\(j\) used to refer to individual visits. The mapping \(\Gamma: \mathbb{V} \rightarrow \mathbb{B}\) is used to map a visit
index, \(i\), to a bus index, \(b\). The notation \(\Gamma_i\) is used as a shorthand to refer to the bus index \(b\) for visit
\(i\).

The actual physical dimensions of the BEB are ignored and it is assumed that each BEB will be assigned to charge at a
particular charger. Because of this assumption, the PAP spatial variable, \(l_i\), may be removed and \(v_i\) is made to be
an integer corresponding to which queue visit \(i\) will be using, \(v_i \in \mathbb{Q}\). That is, the queue position is now
discretized over \(n_Q\) chargers where a BEB occupies single charge queue. Thus, when \(\psi_{ij} = 1\), vehicle \(j\) is placed
in a charging queue with a larger index than vehicle \(i\), \(v_j > v_i\). The charger length \(L\) is likewise replaced with
\(n_Q\). Note that \(n_Q = n_B + n_C\), where \(n_B\) is the number of buses and \(n_C\) is the number of chargers. The
rationale for adding additional idle queues is to allow BEBs to be ``set aside'' if no additional charge is required.
Adding one idle queue for each BEB ensures that the constraints will be satisfied if multiple buses sharing overlapping
times at the station are placed in idle queues. This method will be applied when defining the parameters in
\autoref{sec:example}. The modified queuing constraints can be written as shown in \autoref{eq:packconstrs}.

\begin{subequations}
\label{eq:packconstrs}
\begin{align}
    v_i - v_j - (\psi_{ij} - 1)n_Q \geq 1 \label{subeq:space} \\ d_i \leq \tau_i \label{subeq:valid_depart} \\ s_i \geq
    0 \label{subeq:pos_charge} \\ v_i \in \mathbb{Q} \label{subeq:vspace}
\end{align}
\end{subequations}

The constraint in \autoref{subeq:space} is nearly identical to \autoref{subeq:bapspace}, but rather than viewing the
charger as a continuous strip of length \(S\), it is discretized into \(n_Q\) queues each with a width of unit length one. A
BEB is also assigned a unit length of one which is reflected in \autoref{subeq:space} by \(\cdot \geq 1\).
\autoref{subeq:valid_depart} ensures that the time the BEB is detached from the charger, \(d_i\), is before its departure
time, \(\tau_i\) seconds. Note the introduction of the new variable \(\tau_i\) exists to allow the final charge time to be
independent a similar manner that the inital charge time need not coincide with the arrival time, \(a_i \le u_i \le d_i \le
\tau_i\). \autoref{subeq:vspace} defines the of integers that \(v_i\) that represent the \(n_Q\) chargers.

\subsection{Battery Charge Dynamic Constraints}
\label{sec:batt_dynamics}
Battery dynamic constraints are now to be introduced. Two constraints are enforced on the SOC for each BEB: the SOC must
always remain above a specified percentage to guarantee sufficient charge to execute their respective routes and each
bus must end the day with an SOC above a specified threshold, preparatory for the next day.

The SOC upon arrival for visit \(i\) is denoted as \(\eta_i\) kWh. Because the SOC for a visit \(i\) is dependent on its previous
visits, the mapping \(\Upsilon: \mathbb{V} \rightarrow \mathbb{V} \bigcup \{\varnothing\}\) is used to determine the next visit that corresponds
to the same bus, with \(\Upsilon_i\) being shorthand notation. Thus, \(\Gamma_j\) and \(\Gamma_{\Upsilon_i}\), for \(\Upsilon_i = j\), would both map to the
same bus index as long as \(\Upsilon_i\) is not the null element, \(\varnothing\). The null element is reserved for BEBs that have
no future visits.

To drive time spent on the charger, $s_i$, as well as define initial, final, and intermediate bus charges for each visit
$i$, the sets for initial and final visits must be defined. Let the mapping of the first visit by each bus be denoted as
$\Gamma^0 : \mathbb{B} \rightarrow \mathbb{V}$. The resulting value of the mapping $\Gamma^0$ represents the index for
the first visit of bus $b$. Similarly, let $\Gamma^f : \mathbb{B} \rightarrow \mathbb{V}$ maps the indices for the final
visits for each bus $b \in \mathbb{B}$. Let the storthand for each mapping be denoted as $\Gamma^0_b$ and $\Gamma^f_b$,
respectively. The initial and final bus charge percentages, $\alpha$ and $\beta$, can then be represented by the
constraint equations $\eta_{\Gamma^0_b} = \alpha_b \kappa_{b}$ and $\eta_{\Gamma^f_b} = \beta_b \kappa_{b}$,
respectively. The intermediate charges must be determined during runtime.

It is assumed that the charge received is proportional to the time spent charging. The rate for charger \(q\) is denoted
as \(r_q\) kW. Note that a value of \(r_q = 0\) corresponds to a queue where no charging occurs. A bus in such a queue is
simply waiting at the station for the departure time. The queue indices are ordered such that the first \(n_B\) queues
have \(r_q = 0\) to allow an arbitrary number of buses to sit idle at any given moment in time. The next \(n_C\) queues are
reserved for the slow and fast chargers. The amount of discharge between visits \(i\) and \(\Upsilon_i\), the next visit of the
same bus, is denoted as \(\Delta_i\) kWh. If visit \(i\) occurred at charger \(q\), the SOC of the BEB's next arrival, \(\Upsilon_i\), would
be \(\eta_{\Upsilon_i} = \eta_i + s_i r_q - \Delta_i\).

The binary decision variable \(w_{iq} \in \{0,1\}\) is introduced to indicate the active charger for visit \(i\) in vector
form. The form of the SOC for the next visit, \(\Upsilon_i\), can be written using the following constraints.

\begin{subequations}
    \label{subeq:pre_next_charge}
\begin{align}
    \eta_{\Upsilon_i} = \eta_i + \sum_{q=1}^{n_Q} s_i w_{iq} r_q - \Delta_i \\
    \sum_{q=1}^{n_Q} w_{iq} = 1                           \\
    w_{iq} \in \{0,1\}.
\end{align}
\end{subequations}

The choice of queue for visit \(i\), becomes a slack variable and is defined in terms of \(w_{iq}\) as

\begin{equation}
    v_i = \sum_{q=1}^{n_Q} qw_{iq}.
\end{equation}

Maximum and minimum values for the charges are included to ensure that the battery is not overcharged and to guarantee
sufficient charge for subsequent visits. The upper and lower battery charge bounds for bus \(b\) are \(\kappa_b\) and \(\nu_b \kappa_b\),
respectively , where \(\kappa_b\) is the battery capacity and \(\nu_b\) is a percent value. The upper and lower bounds for the
current SOC are written as follows.

\begin{subequations}
    \label{subeq:pre_min_max}
\begin{align}
    \eta_i + \sum_{q=1}^{n_Q} s_i w_{iq} r_q \leq \kappa_{\Gamma_i} \label{eq:maxcharge}\\
    \eta_i \geq \nu_{\Gamma_i} \kappa_{\Gamma_i} \label{eq:mincharge}
\end{align}
\end{subequations}

\autoref{eq:maxcharge} ensures that the BEB SOC does not exceed the battery capacity, and \autoref{eq:mincharge}
enforces that the inital SOC for each visit is above the threshold of \(\nu_{\Gamma_i}\kappa_{\Gamma_i}\). Note that the term \(s_i w_{iq}\)
is a bilinear term. A standard way of linearizing a bilinear term that contains an integer variable is by introducing a
slack variable with an either/or constraint \cite{chen-2010-applied,rodriguez-2013-compar-asses}. Allowing the slack
variable \(g_{iq}\) seconds to be equal to \(s_i w_{iq}\), \(g_{iq}\) can be defined as

\begin{equation}
    \label{eq:giq_cases}
    g_{iq} =
    \begin{cases}
        s_i & w_{iq} = 1 \\
        0 & w_{iq} = 0
    \end{cases}.
\end{equation}

\autoref{eq:giq_cases} can be expressed as a mixed integer constraint using big-M notation with the following four
constraints.

\begin{subequations}
    \label{eq:slack_gain}
\begin{align}
    s_i - (1 - w_{iq})M \leq g_{iq}  \label{subeq:repgpgret} \\
    s_i \geq g_{iq}                 \label{subeq:repgples} \\
    Mw_{iq} \geq g_{iq}              \label{subeq:repgwgret} \\
    0 \leq g_{iq}                   \label{subeq:repgwles}
\end{align}
\end{subequations}

\noindent where \(M\) is a large unitless value. If \(w_{iq} = 1\) then \autoref{subeq:repgpgret} and
\autoref{subeq:repgples} become \(s_i \leq g_{iq}\) and \(s_i \geq g_{iq}\), forcing \(s_i = g_{iq}\) with \autoref{subeq:repgwgret}
being inactive. If \(w_{iq} = 0\), \autoref{subeq:repgpgret} is inactive and \autoref{subeq:repgwgret} and
\autoref{subeq:repgwles} force \(g_{iq} = 0\).

\subsection{The BEB Charging Problem}
\label{sec:BEB_MILP}
The goal of the MILP is to utilize chargers as little as possible to reduce energy costs with fast charging being
penalized more to avoid the adverse effects of fast charging on battery health as well as the
larger usage cost. Thus, an assignment cost \(m_q\) and usage cost \(\epsilon_q\) are associated with each charger, \(q\).
These unitless weights can be adjusted based on charger type or time of day that the visit
occurs. The assignment term takes the form \(w_{iq}m_q\), and the usage term takes the form \(g_{iq} \epsilon_q\). The
resulting BEB charging problem is defined in \autoref{eq:objective}.

\begin{equation}
\label{eq:objective}
	\min \sum_{i=1}^N \sum_{q=1}^{n_Q} \Big( w_{iq} m_q + g_{iq} \epsilon_q \Big) \\
\end{equation}

Subject to the constraints

\begin{multicols}{2}
\begin{subequations}
                                                     \label{eq:dynconstrs}
\begin{equation}
    u_j - u_i - s_i - (\sigma_{ij} - 1)T \geq 0              \label{subeq:m_time}         \\
\end{equation}
\begin{equation}
    v_j - v_i - (\psi_{ij} - 1)n_Q \geq 1                  \label{subeq:m_space}        \\
\end{equation}
\begin{equation}
    \sigma_{ij} + \sigma_{ji} + \psi_{ij} + \psi_{ji} \geq 1            \label{subeq:m_valid_pos}    \\
\end{equation}
\begin{equation}
    \sigma_{ij} + \sigma_{ji} \leq 1                              \label{subeq:m_sigma}        \\
\end{equation}
\begin{equation}
    \psi_{ij} + \psi_{ji} \leq 1                              \label{subeq:m_delta}        \\
\end{equation}
\begin{equation}
    s_i + u_i = d_i                                  \label{subeq:m_detach}       \\
\end{equation}
\begin{equation}
    \eta_{\Gamma^0_b} = \alpha_{\Gamma_i} \kappa_{\Gamma_i}                         \label{subeq:init_charge}    \\
\end{equation}
\begin{equation}
    a_i \leq u_i \leq (T - s_i)                            \label{subeq:m_valid_starts} \\
\end{equation}
\begin{equation}
    d_i \leq \tau_i                                        \label{subeq:m_valid_depart} \\
\end{equation}
\begin{equation}
    \eta_i + \sum_{q=1}^{n_Q} g_{iq} r_q - \Delta_i = \eta_{\gamma_i}   \label{subeq:next_charge}    \\
\end{equation}
\begin{equation}
    \eta_i + \sum_{q=1}^{n_Q} g_{iq} r_q - \Delta_i \geq \nu_{\Gamma_i} \kappa_{\Gamma_i} \label{subeq:min_charge}     \\
\end{equation}
\begin{equation}
    \eta_i + \sum_{q=1}^{n_Q} g_{iq} r_q \leq \kappa_{\Gamma_i}         \label{subeq:max_charge}     \\
\end{equation}
\begin{equation}
    \eta_{\Gamma^f_b} \geq \beta_{\Gamma_f} \kappa_{\Gamma_f}                   \label{subeq:final_charge}   \\
\end{equation}
\begin{equation}
    s_i - (1 - w_{iq})M \leq g_{iq}                     \label{subeq:gpgret}         \\
\end{equation}
\begin{equation}
    s_i \geq g_{iq}                                     \label{subeq:gples}          \\
\end{equation}
\begin{equation}
    Mw_{iq} \geq g_{iq}                                 \label{subeq:gwgret}         \\
\end{equation}
\begin{equation}
    0 \leq g_{iq}                                       \label{subeq:gwles}          \\
\end{equation}
\begin{equation}
    v_i = \sum_{q=1}^{n_Q} qw_{iq}                      \label{subeq:wmax}           \\
\end{equation}
\begin{equation}
    \sum_{q=1}^{n_Q} w_{iq} = 1                         \label{subeq:wone}           \\
\end{equation}
\begin{equation}
   w_{iq}, \sigma_{ij}, \psi_{ij} \in \{0,1\}\;            \label{subeq:binaryspace}        \\
\end{equation}
\begin{equation}
    v_i, q_i \in  \mathbb{Q}                                         \label{subeq:Qspace}        \\
\end{equation}
\begin{equation}
    i \in \mathbb{V}                                   \label{subeq:Ispace}         \\
\end{equation}
\end{subequations}
\end{multicols}

\autoref{subeq:m_time}-\autoref{subeq:m_valid_depart} are reiterations of the queuing constraints in
\autoref{eq:bapconstrs}. \autoref{subeq:init_charge}-\autoref{subeq:final_charge} provide the battery charge
constraints. \autoref{subeq:gpgret}-\autoref{subeq:gwles} define the charge gain of every visit/queue pairing. The last
constraints \autoref{subeq:binaryspace}-\autoref{subeq:Ispace} define the sets of valid values for each variable.
\section{Example}
\label{sec:example}
An example will now be presented to demonstrate the utility of the developed MILP charge scheduling technique. A
description of the scenario is first presented followed by a description of an alternative heuristic-based planning
strategy called Qin-Modified which is used as a comparison to the MILP PAP. Results are then presented,
analyzed and discussed for each of the planning strategies.

\subsection{BEB Scenario}
\label{beb-scenario}
To display the capabilities of the model, an example scenario is presented. The scenario was run over a time horizon of
\(T=24\) hours, utilizing \(n_B = \A\) buses with \(n_V = \N\) visits divided between the \(n_B\) buses. As stated before, the
route times are sampled from a set of routes from the UTA. Each bus has a battery capacity of \(\kappa_b =\) \batsize kWh that
is required to stay above an SOC of \(\nu_b =\) \mincharge (\fpeval{\batsize * \minchargeD} kWh) \(\forall b \in
\mathbb{B}\) to ensure that each BEB can complete its next route in addition to maintaining battery health. Each bus is
assumed to begin the working day with an SOC of \(\alpha_b =\) \fpeval{\acharge*100}\% \(\forall i \in \mathbb{V}\)
(\fpeval{\acharge * \batsize} kWh). Additionally, each bus is required to end the day with a minimum SOC of \(\beta_b =\)
\fpeval{\bcharge * 100}\% (\fpeval{\bcharge * \batsize} kWh) \(\forall i \in \mathbb{V}\). This assumes that overnight
charging can account for the deficient 20\% SOC. Each bus is assumed to discharge at a rate of \(\zeta_b =\) 30 kW. Note that
many factors play a role in the rate of discharge; however, for the sake of simplicity since the discharge calculation
is out of the scope of this work, an average rate is used. A total of \(n_C = \fpeval{\fast + \slow}\) chargers are
utilized where \slow of the chargers are slow charging (\slows kW) and \fast are fast charging (\fasts kW). The
technique to minimize the total charger count will now be employed.

To encourage the MILP PAP problem to utilize the fewest number of chargers, the value of \(m_q\) in the objective
function, \autoref{eq:objective}, is \(\forall q \in \{1,2,..., n_B \}; m_q = 0\) and \(\forall q \in \{n_B + 1, n_B + 2,..., n_B + n_C \};
m_q = 1000q\). The charge duration scalar, \(\epsilon_q\), is defined as \(\epsilon_q = r_q\) to create a consumption cost term,
\(g_{iq}\epsilon_q\) kWh. By consumption cost, it is meant that the total energy consumed by the charge schedule will be
accounted for in the objective function. This method encourages the model to minimize active charger times, particularly
for the fast chargers.

Another heuristic-based optimization strategy, referred to as Qin-Modified, is also employed as a means of comparison
with the results of the MILP PAP. The Qin-Modified strategy is based on the threshold strategy of
\cite{qin-2016-numer-analy}. The strategy has been modified slightly to accommodate the case of multiple charger types
without an exhaustive search for the best charger type. The heuristic is based on a set of rules that revolve around the
initial SOC of the bus visit \(i\). There are three different thresholds, low (85\%), medium (90\%), and
high (95\%). Buses below the low threshold (\(\text{SOC} \le 85\%\)) are prioritized to fast chargers and then are allowed to
utilize slow chargers if no fast chargers are available. Buses between the low and medium threshold (\(85\% < \text{SOC}
\le 90\%\)) prioritize slow chargers first and utilizes fast chargers only if no slow chargers are available. Buses above
the medium threshold and below the high (\(90\% < \text{SOC} \le 95\%\)) will only be assigned to slow chargers. Buses above
the high threshold (\(\text{SOC} > 95\%\)) will not be placed in a charging queue. Once a bus has been assigned to a
charger, it remains on the charger for the duration of the time it is at the station, or it reaches an SOC of 95\%
charge, whichever comes first.

The total number of constraints resulted in \contvars continuous and \intvars integer/binary constraints. The
optimization was performed using the Gurobi MILP solver \cite{gurobi-2021-gurob-optim} on a machine equipped with an
AMD Ryzen 9 5900X 12 - Processor (24 core) at 4.95GHz. The solver was allowed to run for \timeran seconds.

\subsection{Results}
\label{results}
The schedule generated by the Qin-Modified strategy and the MILP PAP is shown in \autoref{subfig:qin-schedule} and
\autoref{subfig:milp-schedule}, respectively. The x-axis represents the time in hours. The y-axis represents the
assigned charging queue. Rows between 0 and \fpeval{\slow - 1} are active times for slow chargers, and rows in the range
of \fpeval{\slow} and \fpeval{\fast + \slow - 1} are active times for fast chargers. The unique color/symbol-styled
vertices represent the starting charge time for a bus \(b\) with the line to the vertical tick signifying the region of
time the charger is active. The lines connecting points represents the charge sequence for each BEB.

The first observation is in the choice of preferred chargers between the Qin-Modified and MILP scheduler. Looking at
\autoref{subfig:slow-charger-usage} and \autoref{subfig:fast-charger-usage}, the Qin-Modified schedule uses at most four
fast chargers and three slow at the same time, whereas the MILP schedule uses at most one fast charger and six slow at
the same time. Both the Qin-Modified and MILP schedule used the fast chargers in short bursts (\textasciitilde{}0.2-0.5 hours). The main
difference lies in the utilization strategy of the slow chargers. The Qin-Modified, for the most part, opted for shorter
bursts for the slow chargers (\textasciitilde{}0.3-0.7 hours), most heavily placed on the first slow charger. The MILP
utilized the slow chargers in short bursts; however, the solver was able to recognize moments where a BEB being placed
in a slow charging queue for a longer duration was more cost-effective (with respect to the objective function) than
placing the BEB in a fast charging queue. Although one of the MILP's objectives is to minimize the amount of chargers
used, the Qin-Modified schedule ended up using fewer chargers than the MILP. Note the MILP schedule
packed the first queue for the fast and slow chargers more effectively than the Qin schedule. Although both schedules
generated are valid, no comparison of the quality of the schedule can be made directly from
\autoref{subfig:milp-schedule} and \autoref{subfig:qin-schedule}.

\autoref{subfig:qin-charge} and \autoref{subfig:milp-charge} depicts the SOC for every bus over the time horizon for the
Qin and MILP schedules, respectively. Every vehicle begins with an SOC of \(\alpha_b =\) 90\%, finishes with an SOC of \(\beta_b =\)
\fpeval{\bcharge *100}\% in the MILP PAP schedule, and never goes below \mincharge in the intermediate arrivals as stated
in constraint \autoref{eq:dynconstrs}. There is no guarantee for this in the Qin-Modified strategy which can be seen by
some intermediate charges reaching an SOC of 0\% as well as the distribution of final charges, the minimum being 0\% and
the maximum \fpeval{trunc(\fpeval{368 / \batsize * 100}, 3)}\%. The only sense of guarantee that the
Qin-Modified supplies is its predictability within the intermediate visits due to its heuristic nature (i.e. if the BEB
charge is within the low threshold, a fast charger will be prioritized); whereas MILP places a bus in the queue that
``makes sense'' in respect to the larger picture. The MILP PAP does not have an obvious sense of decision-making due to
weighted objective function that is affected by the accumulation of decisions made prior.

\begin{figure}[htpb]
\centering
    \includegraphics{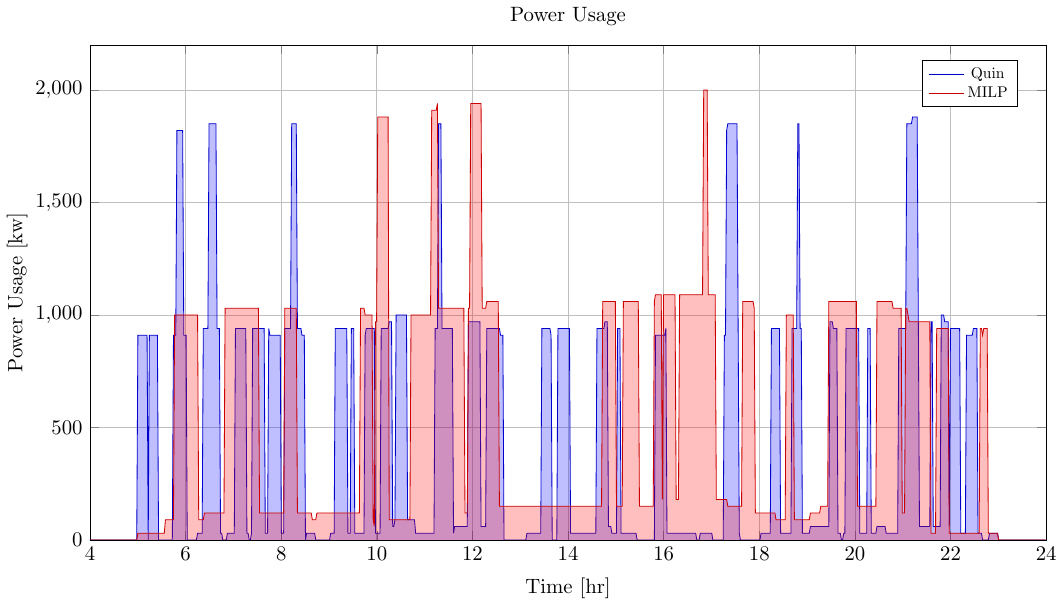}
    \caption{Amount of power consumed by Qin-Modified and MILP schedule over the time horizon.}
    \label{fig:power-usage}
\end{figure}

\begin{figure}[htpb]
\centering
    \includegraphics{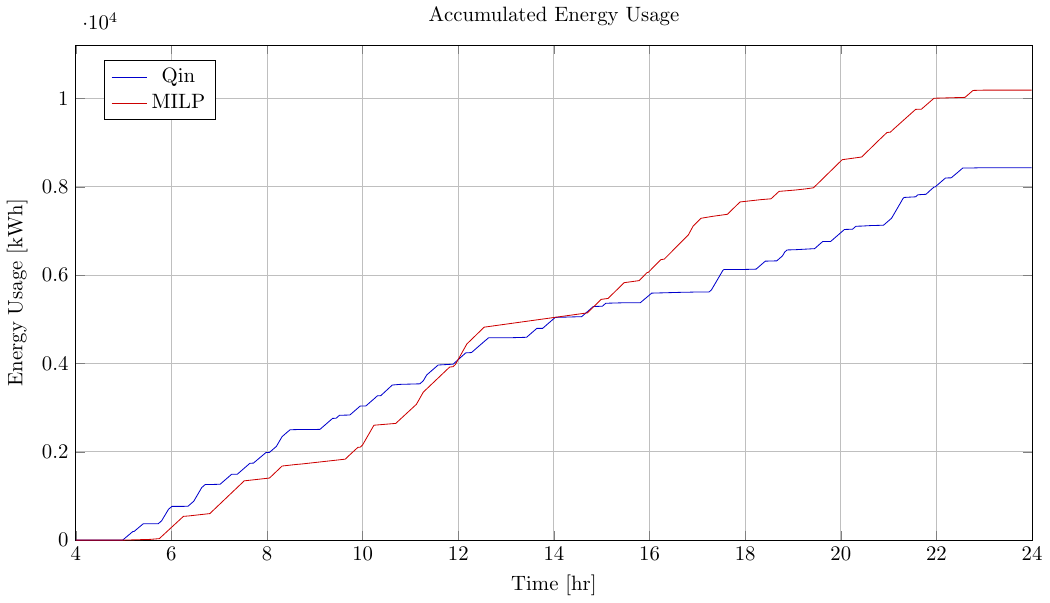}
    \caption{Total accumulated energy consumed by the Qin-Modified and MILP schedule throughout the time horizon.}
    \label{fig:energy-usage}
\end{figure}

Another important measure for the chargers is to compare the amount of power and energy consumed.
\autoref{fig:power-usage} depicts the power consumption throughout the time horizon. It can be seen that the
Qin-Modified power consumption is steadily less or the same as the MILP schedule. This can be accounted for by the
MILP's constraints to keep the bus SOC above \mincharge and to reach a final SOC of \fpeval{\bcharge *100}\% at the end
of the working day. Along a similar vein, the accumulated energy consumed is shown in \autoref{fig:energy-usage}. The
MILP schedule is more efficient up until about the eleventh hour. Again, this can be accounted for by the fact the MILP
is accommodating the extra constraints. Due to these constraints the MILP PAP consumes about \(0.1\cdot10^4\) kWh more than
the Qin-Modified. The overlap of the MILP PAP can be accounted for by referencing \autoref{subfig:fast-charger-usage}
and \autoref{subfig:slow-charger-usage}. Between the fifth and tenth hour, the MILP schedule heavily uses slow chargers
increasing the rate at which power is being consumed. Afterwards, the MILP schedule at a minimum continues to use the
same amount of chargers as the Qin Schedule. Again, due to the added constraints, the MILP schedule must utilize more
resources to keep within the specified bounds.
\section{Conclusion}
\label{sec:conclusion}
This work developed a MILP scheduling framework that optimally assigns fast and slow chargers to a BEB bus fleet
assuming a constant schedule. The BAP was briefly introduced followed by a description and formulation of the PAP. The
PAP was modified to allow charge time to be variable. Because the modified PAP no longer requires the charge time to be
given, linear battery dynamics were introduced to propagate each BEBs' SOC. Additional constraints were also introduced
to provide upper and lower limits for the battery dynamics.

An example for the MILP PAP formulation was then presented and compared to a heuristic-based schedule, referred to as
Qin-Modified. The MILP PAP optimization was run for \timeran seconds to a non-optimal solution. The Qin-Modified and
MILP schedule utilized a similar amount of fast chargers; however, the MILP schedule more readily used the slow chargers
to its advantage when the objective function saw fit. More importantly, the MILP PAP schedule utilized approximately
\(0.1\cdot10^4\) kWh more than the Qin-Modified; however, the charges for the MILP schedule remained above the constrained
minimum SOC of \mincharge, and charged all the buses to \fpeval{\bcharge *100}\% at the end of the working day. The
Qin-Modified schedule, on the other hand, allowed the SOC of certain BEBs to drop to 0\%.

Further fields of interest are to utilize the formulation (\autoref{eq:objective} and \autoref{eq:dynconstrs}) with
nonlinear battery dynamics, calculation and utilization of the demand and consumption cost in the objective function,
and utilizing this formulation in a metaheuristic solver. Furthermore, ``fuzzifying'' the initial and final charge times
is of interest to allow flexibility in the arrival and departure times.

\bibliographystyle{Frontiers-Harvard}
\bibliography{/home/alex/Documents/citation-database/lit-ref,/home/alex/Documents/citation-database/lib-ref}

\clearpage

\section{Figure Captions}
\label{sec:orgfdd54a7}

\begin{subfigures}
    \begin{figure}[htpb]
    \centering
        \includegraphics{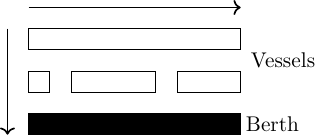}
        \caption{Example of berth allocation. Vessels are docked in berth locations (horizontal) and are queued over
          time (vertical). The vertical arrow represents the movement direction of queued vessels and the horizontal
          arrow represents the direction of departure.}
        \label{subfig:bapexample}
    \end{figure}
    \hfill

    \begin{figure}[htpb]
    \centering
        \includegraphics{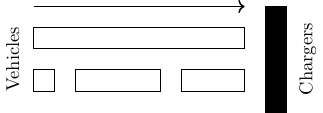}
        \caption{Example of position allocation. Vehicles are placed in queues to be charged and move in the direction
          indicated by the arrow.}
        \label{subfig:papexample}
    \end{figure}
\end{subfigures}

\begin{table}[!htpb]
  \caption{Notation used throughout the paper. Units are provided when available.}
  \label{tab:variables}
  \centering
  \begin{tabularx}{\textwidth}{l l l}
    \toprule \textbf{Variable} & \textbf{Units} & \textbf{Description}                                                                      \\
    \toprule \multicolumn{3}{l}{Input values}                                                                                               \\
    \hline $n_B$ & & Number of buses                                                                                                        \\
    $M$           &       & An arbitrarily large number                                                                                     \\
    $n_V$         &       & Number of total visits                                                                                          \\
    $n_Q$         &       & Number of queues                                                                                                \\
    $n_C$         &       & Number of chargers                                                                                              \\
    $\mathbb{V}$  &       & Set of visit indices, $\mathbb{V} = \{1, ..., n_V\}$                                                            \\
    $\mathbb{B}$  &       & Set of bus indices, $\mathbb{B} = \{1, ..., n_B\}$                                                              \\
    $\mathbb{Q}$           &       & Set of queue indices, $\mathbb{Q} = \{1, ..., n_Q\}$                                                                     \\
    $i,j$         &       & Indices used to refer to visits                                                                                 \\
    $b$           &       & Index used to refer to a bus                                                                                    \\
    $q$           &       & Index used to refer to a queue                                                                                  \\
    \hline \multicolumn{3}{l}{Problem definition parameters}                                                                                \\
    \hline $\Gamma$    &       & $\Gamma: \mathbb{V} \rightarrow \mathbb{B}$ with $\Gamma_i$ used as a shorthand to denote the bus $b$ for visit $i$                 \\
    $\alpha_b$         & $\%$  & Initial charge percentage time for bus $b$                                                                      \\
    $\beta_b$         & $\%$  & Final charge percentage for bus $b$ at the end of the time horizon                                              \\
    $\epsilon_q$         &       & Cost of using charger $q$ per unit time                                                                         \\
    $\Upsilon$           &       & $\Upsilon: \mathbb{V} \rightarrow \mathbb{V}$ mapping a visit to the next visit by the same bus with $\Upsilon_i$ being the shorthand.  \\
    $\kappa_b$         & kWh   & Battery capacity for bus $b$                                                                                    \\
    $\Delta_i$         & kWh   & Discharge of visit over route $i$                                                                               \\
    $\nu_b$         & $\%$  & Minimum charge allowed for bus $b$                                                                              \\
    $\tau_i$         & s     & Time visit $i$ must depart the station                                                                          \\
    $\zeta_b$         & kW    & Discharge rate for bus $b$                                                                                      \\
    $a_i$         & s     & Arrival time of visit $i$                                                                                       \\
    $i_0$         &       & Indices associated with the initial arrival
    for every bus in $\mathbb{B}$                                                                                                           \\
    $i_f$         &       & Indices associated with the final arrival for every bus in $\mathbb{B}$                                         \\
    $m_q$         &       & Cost of a visit being assigned to charger $q$                                                                   \\
    $r_q$         & kW    & Charge rate of charger $q$ per unit time                                                                        \\
    \hline \multicolumn{3}{l}{Decision Variables}                                                                                           \\
    \hline $\psi_{ij}$          &                         & Binary variable determining spatial ordering of vehicles $i$ and $j$               \\
    $\eta_i$     & kWh  & Initial charge for visit $i$                                                                                         \\
    $\sigma_{ij}$  &      & Binary variable determining temporal ordering of vehicles $i$ and $j$                                                \\
    $d_i$     & s    & Ending charge time for visit $i$                                                                                     \\
    $g_{iq}$  & s    & Linearization term, represents the multiplication of $s_i w_{iq}$                                                     \\
    $s_i$     & s    & Amount of time spent on charger for visit $i$                                                                        \\
    $u_i$     & s    & Starting charge time of visit $i$                                                                                    \\
    $v_i$     &      & Assigned queue for visit $i$                                                                                         \\
    $w_{iq}$  &      & Binary assignment variable for visit $i$ to queue $q$                                                                \\
    \bottomrule
  \end{tabularx}
\end{table}

\begin{figure}[htpb]
\centering
    \includegraphics{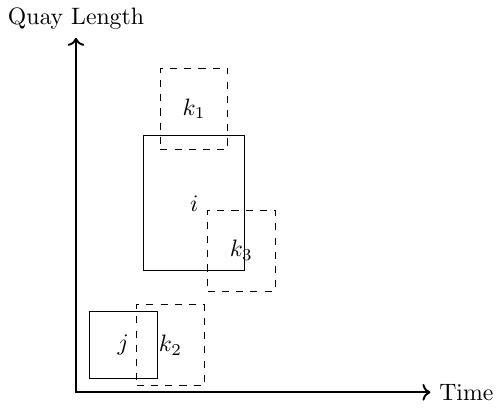}
    \caption{Examples of different methods of overlapping. Space overlap: $v_{k_1} < v_{i} + l_i \therefore \psi_{k_{1}i} = 0$.
             Time overlap $u_{k_1} < u_{j} + s_j \therefore \sigma_{k_{2}j} = 0$. Both space and time overlap $\sigma_{k_{3}i} = 0$ and
             $\psi_{k_{3}j} = 0$.}
    \label{fig:multipleassign}
\end{figure}

\begin{subfigures}
    \begin{figure}[htpb]
    \centering
    \includegraphics{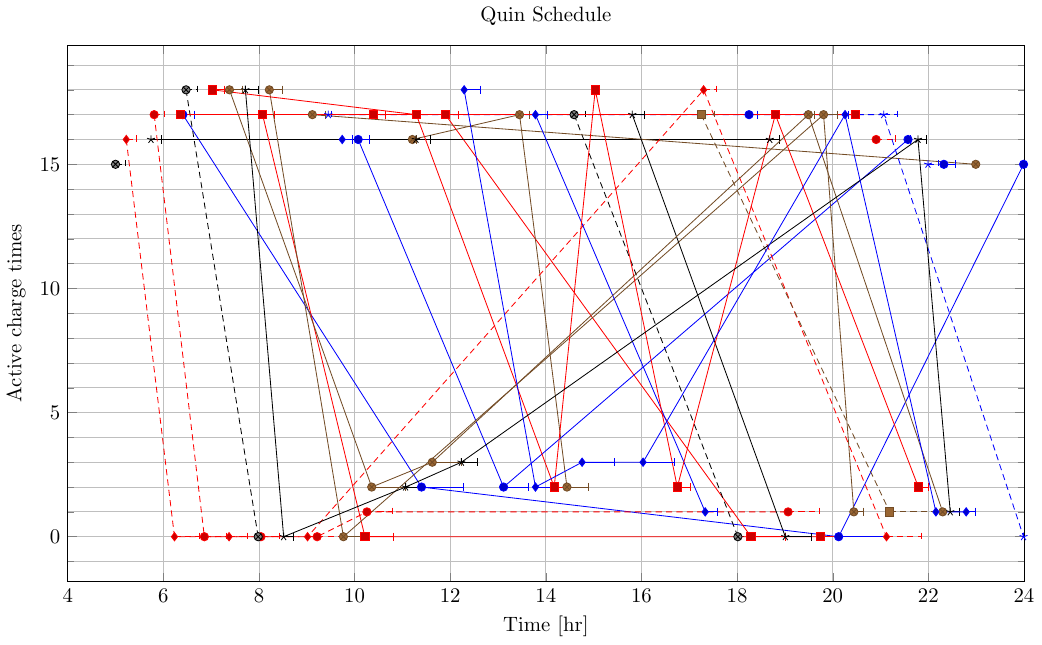}
        \caption{Charging schedule generated by Qin Modified algorithm.}
        \label{subfig:qin-schedule}
    \end{figure}

    \hfill

    \begin{figure}[htpb]
    \centering
        \includegraphics{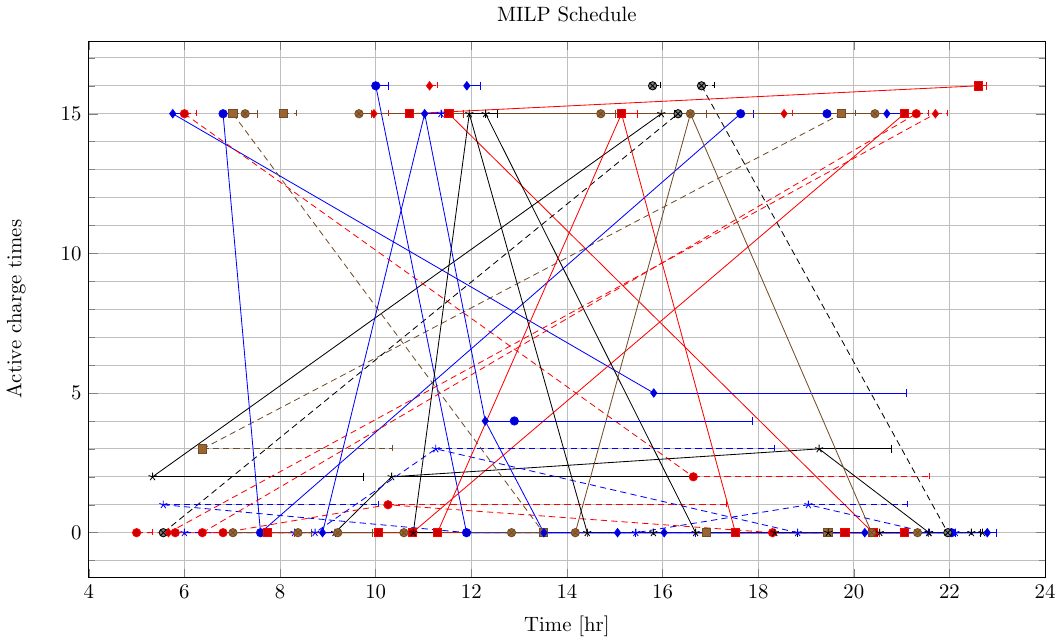}
        \caption{Charging schedule generated by MILP PAP algorithm.}
        \label{subfig:milp-schedule}
    \end{figure}
\end{subfigures}

\begin{subfigures}
    \begin{figure}[htpb]
    \centering
        \includegraphics{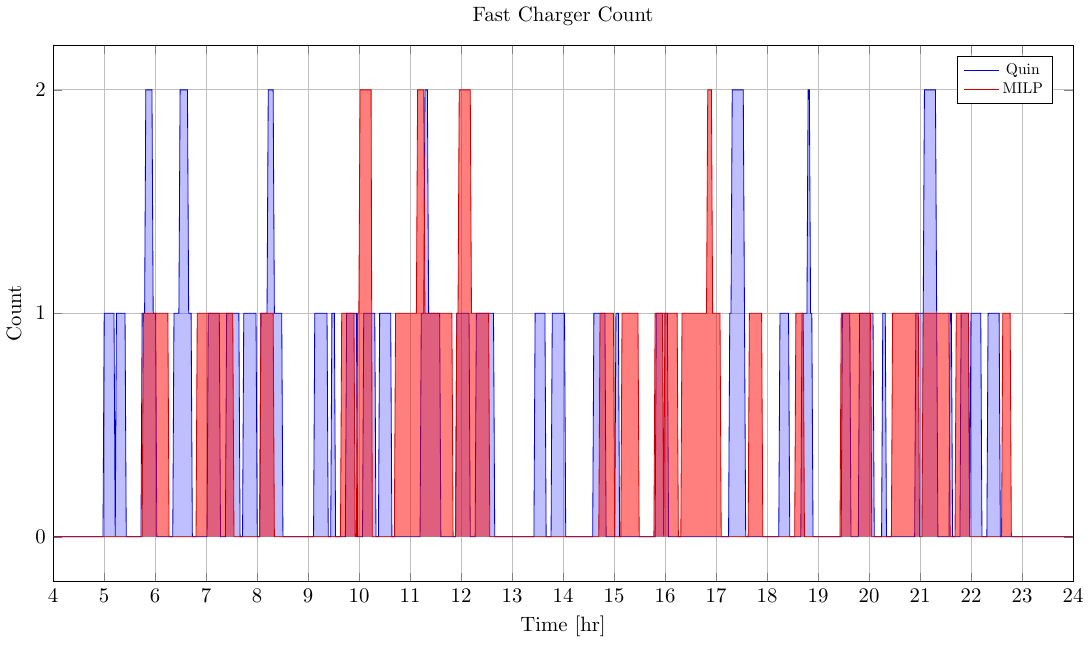}
        \caption{Number of fast chargers for Qin and MILP PAP.}
        \label{subfig:fast-charger-usage}
    \end{figure}

    \hfill

    \begin{figure}[!ht]
    \centering
        \includegraphics{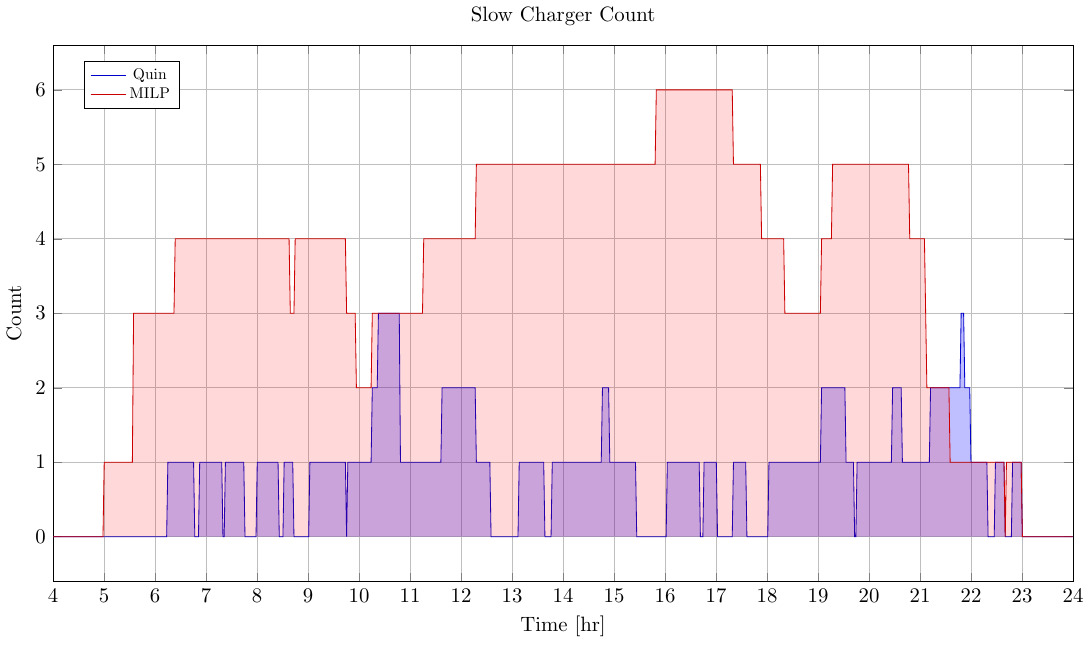}
        \caption{Number of slow chargers for Qin and MILP PAP.}
        \label{subfig:slow-charger-usage}
    \end{figure}
\end{subfigures}

\begin{subfigures}
    \begin{figure}[htpb]
    \centering
        \includegraphics{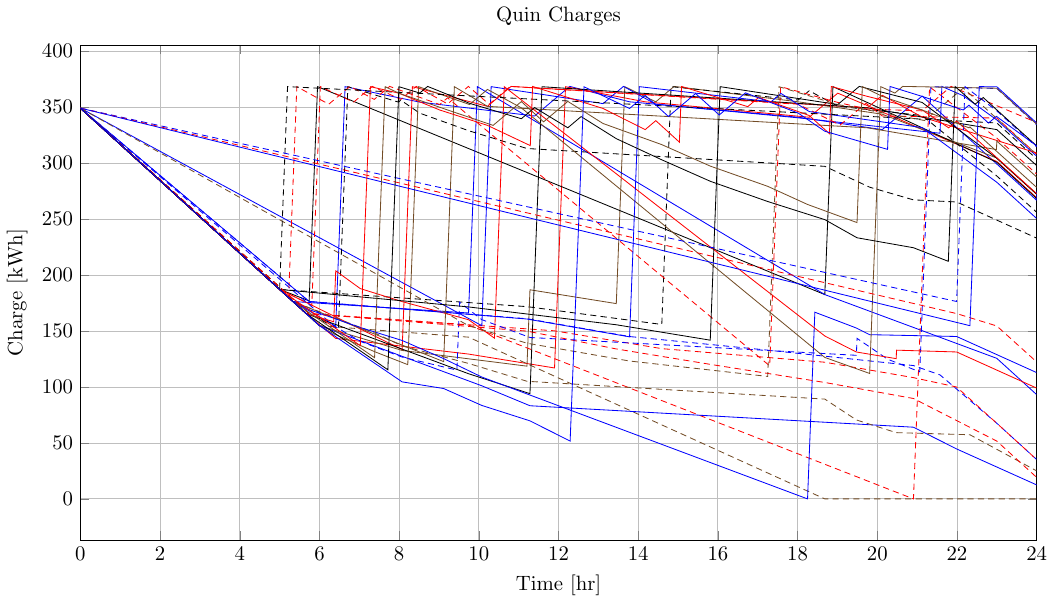}
        \caption{Bus charges for the Qin Modified charging schedule. The charging scheme of the Qin charger is more predictable during the working day.}
        \label{subfig:qin-charge}
    \end{figure}

    \hfill

    \begin{figure}[htpb]
    \centering
        \includegraphics{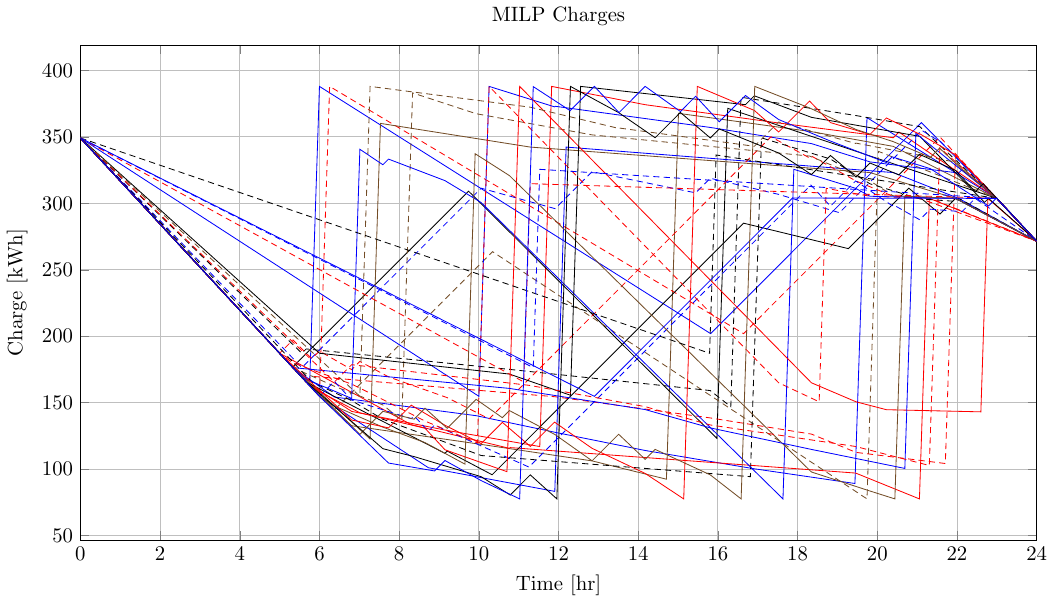}
        \caption{The bus charges for the MILP PAP charging schedule. The MILP model allows for guarantees of minimum/maximum changes during the working day as well as charges at the end of the day.}
        \label{subfig:milp-charge}
    \end{figure}
\end{subfigures}

\end{document}